\documentclass{amsart}
\usepackage{amssymb, bbm, mathrsfs}
\usepackage[all]{xy}
\usepackage{latexsym}

\def\A{{\mathbb A}}
\def\Ii{{\mathbb I}}
\def\sA{{\mathscr A}}
\def\sT{{\mathscr T}}
\def\sS{{\mathscr S}}

\def \im{\mathop{\sf Im}\nolimits}

\def \hom{\mathop{\sf Hom}\nolimits}
\def \ext{\mathop{\sf Ext}\nolimits}
\def\id{\mathop{\sf Id}\nolimits}

\newtheorem{Pro}{Proposition}
\newtheorem{Le}[Pro]{Lemma}
\newtheorem{Th}{Theorem}
\newtheorem{Co}[Pro]{Corollary}

\def\Ab{{\mathbbm{Ab}}}

\begin{document}

\title{Universal coefficient theorem in triangulated categories}

\author{Teimuraz Pirashvili}
\author{Mar\'\i a Julia Redondo}
\thanks{The second author is a researcher from CONICET, Argentina}

\address{Department of Mathematics, University of Leicester, Leicester LE1 7RH, 
United Kingdom.}
\address{Instituto de Matem\'atica, Universidad Nacional del Sur, Av. Alem 1253, (8000) Bah\' \i a Blanca,  Argentina.}


\begin{abstract}
We consider a homology theory $h:\sT \to \sA$ on a triangulated category $\sT$ with values in a graded abelian category $\sA$.
If the functor $h$ reflects isomorphisms, is full and is such that for any object $x$ in $\sA$ there is an object $X$ in $\sT$ with an isomorphism between $h(X)$ and $x$, we prove that $\sA$ is a hereditary abelian category, all idempotents in
$\sT$ split and the kernel of $h$ is a square zero ideal which as a bifunctor on $\sT$ is isomorphic to $\ext^1_{\sA}(h(-)[1], h(-))$.
\end{abstract}

\maketitle

We assume that the reader is familiar with triangulated categories
(see \cite{verdier}, \cite{GM}).   Let us just recall that the triangulated categories were introduced independently by Puppe
\cite{puppe} and by Verdier \cite{verdier}.  Following to Puppe we
do not assume that the octahedral axiom holds.

If $\sT$ is a triangulated
category, the shifting of an object $X\in \sT$ is denoted by
$X[1]$. Assume an abelian category $\sA$ is given, which is
equipped with an auto-equivalence $x\mapsto x[1]$. Objects of
$\sA$ are denoted by the small letters $x,y,z$, etc, while objects
of $\sT$ are denoted by the capital letters $X,Y,Z$, etc. A
\emph{homology theory on $\sT$ with values in $\sA$} is a functor
$h:\sT\to \sA$  such that  $h$ commutes with shifting (up to an
equivalence) and for any distinguished triangle $X\to Y\to Z\to
X[1]$ in $\sT$ the induced sequence $h(X)\to h(Y)\to h(Z)$  is
exact. It follows that then one has the following long exact
sequence
$$\cdots \to h(Z)[-1]\to h(X)\to h(Y)\to h(Z)\to
h(X)[1]\to \cdots$$ In what follows  $\ext^1_\sA(x,y)$ denotes the
equivalence classes of extensions of $x$ by $y$ in the category
$\sA$ and we assume that these classes form a set.\\

In this  paper we prove the following result:

\begin{Th}\label{teorema} Let
$h:\sT\to \sA$ be a homology theory. Assume the following
conditions hold
\begin{itemize}
\item[i)] $h$ reflects isomorphisms, 
\item[ii)] $h$ is full.
\end{itemize}
Then
the ideal
$$\Ii = \{ f \in \hom_{\sT}(X,Y)\mid h(f)=0 \}$$
is a square zero ideal. 
Suppose additionally the following condition holds
\begin{itemize}
\item [iii)] for any short exact sequence $0 \to x \to y \to z \to 0$ in $\sA$ with $x \cong h(X)$ and $z \cong h(Z)$ there is an object $Y \in \sT$ and an isomorphism $h(Y)\cong y$ in $\sA$.
\end{itemize}
Then $\Ii$  is isomorphic as a bifunctor on $\sT$ to 
$$(X,Y)\mapsto \ext^1_{\sA}(h(X)[1],h(Y)).$$
In particular for any $X,Y\in \sT$ one has
the following short exact sequence
$$0\to \ext^1_{\sA}(h(X)[1],h(Y))\to \sT(X,Y)\to \hom_{\sA}(h(X),h(Y))\to 0.$$
Moreover, if we replace condition (iii) by the stronger condition
\begin{itemize}
\item [iv)] for any object $x\in \sA$ there is an object $X\in \sT$ and an
isomorphism $h(X)\cong x$ in $\sA$,
\end{itemize}
then $\sA$ is a hereditary abelian category and all idempotents in
$\sT$ split.
\end{Th}

Thus this is a sort of ''universal coefficient theorem'' in
triangulated categories.\\

Our result is a one step generalization of a well-known result
which claims that if $h$ is an equivalence of categories then
$\sA$ is \emph{semi-simple} meaning that $\ext^1_\sA=0$ (see for
example \cite[p. 250]{GM}).
As was pointed out by J. Daniel Christensen our theorem generalizes Theorem 1.2 and Theorem 1.3 of \cite{phantom} on phantom maps. Indeed let $\sS$  be the homotopy category of spectra or, more generally, a  triangulated category  satisfying axioms 2.1 of \cite{phantom} and let $\sA$ be the category
of additive functors from finite objects of $\sS$ to the category of abelian groups. The category $\sA$
has a shifting, which is given by $(F[1])(X)=F( X[1])$, $F\in \sA$. Moreover let $h:\sS\to \sA$ be a functor given by $h(X)=\pi_0(X\wedge (-))$. Then
$h$ is a homology theory for which the assertions i)-iii) hold and $\Ii(X,Y)$ consists of phantom maps from $X$ to $Y$. Hence by the first part of theorem we obtain the familiar properties of phantom maps. \\

Before we give a proof of the Theorem, let us explain  notations
involved on it.  The functor $h$ reflects isomorphisms, this means
that $f\in \hom_{\sT}(X,Y)$ is an isomorphism provided $h(f)$ is
an isomorphism in $\sA$. This holds if and only if $X=0$ as soon as $h(X)=0$.
Moreover $h$ is full, this means that the homomorphism
$\sT(X,Y)\to \hom_{\sA}(h(X),h(Y))$  given by $f\mapsto h(f)$ is
surjective for all  $X,Y\in \sT$. Furthermore an abelian category
$\sA$ is \emph{hereditary} provided for any two-fold extension
\begin{equation}\label{2gaf}
\xymatrix{0\ar[r]&u\ar[r]^{\hat{\alpha}}&v\ar[r]^{\hat{\beta}}&w \ar[r]^{\hat{\gamma}}
& x\ar[r]&0}
\end{equation} there exists a commutative diagram with exact
rows
$$\xymatrix{&0\ar[r]&v\ar[r]\ar[d]_\id &z\ar[r]\ar[d]&x \ar[r]\ar[d]_\id&0\\
0\ar[r]&u\ar[r] &v\ar[r]&w \ar[r]& x\ar[r]&0}$$ This exactly means
that $\ext^2_\sA=0$, where $\ext$ is understood a l\'a Yoneda. Let
us also recall that an ideal $\Ii$ in an additive category $\A$ is
a sub-bifunctor of the bifunctor $\hom_{\A}(-,-):\A^{op}\times {\A}\to
\Ab$. It follows that $\Ii$ is an additive bifunctor. One can form the quotient category $\A/\Ii$ in
an obvious way, which is an additive category. One says that $\Ii^2=0$ provided $gf=0$ as soon as
$f\in \Ii(A,B)$ and $g\in \Ii(B,C)$. In this case the bifunctor
$\Ii:{\A}^{op}\times \A\to \Ab$ factors through the quotient category
$\A/\Ii$ in a unique way.

\medskip

\begin{proof} It is done in several steps.

\noindent {\it First step. The equality $\Ii^2=0$}.
To make notations easier we denote $h(X), h(Y)$ simply by $x,y$,
etc. Moreover, for a morphism  $\alpha:X\to Y$, we let
$\hat{\alpha}:x\to y$ be the morphism $h(\alpha)$. Suppose $\alpha:X\to Y$
and $\beta:Y\to Z$ are morphisms such that $\hat{\alpha}=0$ and
$\hat{\beta}=0$. We have to prove that $\gamma:=\beta\alpha$ is the zero
morphism. By the morphisms axiom there is a diagram of
distinguished triangles
$$\xymatrix{X\ar[r]^\alpha \ar[d]^{\id}& Y \ar[r] \ar[d]^\beta &U\ar[r]\ar@{-->}[d]
& X[1]\ar[d]^\id\\
X\ar[r]^\gamma& Z\ar[r]^\omega  &V\ar[r]^\nu & X[1]}
$$ Apply $h$ to get a commutative diagram with exact rows
$$\xymatrix{0\ar[r]& y\ar[r] \ar[d]^0 &u\ar[r]\ar[d]
& x[1]\ar[d]^\id\ar[r]& 0\\
0\ar[r]& z\ar[r]^{\hat{\omega}}  &v\ar[r]^{\hat{\nu}} & x[1]\ar[r]&0}$$ It
follows that there is a morphism $\hat{\mu}:x[1]\to v$ in $\sA$ such that
$\hat{\nu}\hat{\mu}=\id_{x[1]}$. Thus $(\hat{\omega},\hat{\mu}):z\oplus
x[1]\to v$ is an isomorphism. Since $h$ is full, we can find $\mu:X[1]\to V$
which realizes $\hat{\mu}$, meaning that $h(\mu)=\hat{\mu}$. The morphism
$(\omega, \mu):Z\oplus X[1]\to V$ is an isomorphism, because $h$  reflects
isomorphisms. In particular $\omega$ is a monomorphism and therefore
$\gamma=0$ and first step is done.

For objects $X,Y\in \sT$ we put
$$\Ii(X,Y):=\{\alpha\in \hom_\sA(X,Y)\mid h(\alpha)=0\}.$$
We have just proved that $\Ii^2=0$. In particular $\Ii$ as a
bifunctor factors through the category $\sT/\Ii$. The next step
shows that it indeed factors through the category $\sA$ and a
quite explicit description of this bifunctor is given.

\smallskip
\noindent {\it Second step. Bifunctorial isomorphism $\Ii(X,Y)\cong
\ext^1_{\sA}(h(X)[1],h(Y))$}.
We put as usual $x=h(X)$, $y=h(Y)$, etc. Let
$\alpha:X\to Y$ be an element of $\Ii(X,Y)$. Consider a distinguished
triangle
\begin{equation}\label{150105}
X\buildrel \alpha \over \to Y\buildrel \beta \over \to Z \buildrel \gamma
\over \to X[1].\end{equation} By applying $h$ one obtains the
following short exact sequence
\begin{equation}\label{1501051}
0\to y \buildrel \hat{\beta}\over \to z \buildrel \hat{\gamma} \over
\to x[1]\to 0\end{equation} whose class in $\ext^1_\sA(x[1],y)$ is
independent on the choice of the triangle in (\ref{150105}) and it
is denoted by $\Xi(\alpha)$. In this way one obtains the binatural
transformation $\Xi: \Ii\to \ext^1_{\sA}((-)[1],(-))$. We claim
that $\Xi$ is an isomorphism. Indeed, if $\Xi(\alpha)=0$, then there
exists a section $\hat{\mu}:x[1]\to z$ of $\hat{\gamma}$ in
(\ref{1501051}). Then $(\hat{\beta},\hat{\mu}):y\oplus x[1]\to z$
is an isomorphism. Since $h$ is full, we can find $\mu:X[1]\to Z$
which realizes $\hat{\mu}$. The morphism $(\beta, \mu):Y\oplus
X[1]\to Z$ is an isomorphism, because $h$  reflects isomorphisms.
In particular $\beta$ is a monomorphism and therefore $\alpha=0$. Hence
$\Xi$ is a monomorphism. Let us take any element in
$\ext^1_{\sA}(x[1],y)$, which is represented by a short exact
sequence, say the sequence (\ref{1501051}). Take any realization
$\beta:Y\to Z$ of $\hat{\beta}$.  By Lemma \ref{mzmrelizacia} below we obtain the following distinguished triangle
$$ X \buildrel \alpha\over \to Y \buildrel \beta \over \to Z \buildrel \gamma \over \to
X[1]
$$
containing $\beta$. It follows that $\Xi(\alpha)$
represents our original element  in $\ext^1_{\sA}(x[1],y)$. Hence
$\Xi$ is an isomorphism.

\smallskip
\noindent {\it Third step. $\sA$ is  hereditary}. Let (\ref{2gaf}) be a two-fold extension in $\sA$. We put
$y=\im(\hat{\alpha})$. Thus the exact sequence (\ref{2gaf}) splits in the
following  two short exact sequences
$$ 0 \to u \buildrel {\hat{\alpha}} \over \to v \buildrel {\hat{\mu}} \over \to y \to 0$$
and
$$0\to y \buildrel {\hat{\nu}} \over \to w \buildrel {\hat{\gamma}} \over \to x\to 0$$
with $\hat{\beta}=\hat{\nu}\hat{\mu}$. Using Assumption iii) and
without loss of generality we can assume that $u,v,w,x$ as well as
$\hat{\alpha}$ and $\hat{\gamma}$ have realizations. By Lemma
\ref{mzmrelizacia} below we obtain the following distinguished
triangles
$$U \buildrel {\alpha} \over \to V \buildrel {\mu} \over \to Y \buildrel {\xi} \over \to U[1]
$$
and
$$
Y \buildrel {\nu} \over \to  W \buildrel {\gamma}\over \to X \buildrel {\chi} \over \to Y[1].
$$
Since $\hat{\mu}$ is an epimorphism and $\hat{\nu}$ is a
monomorphism it follows that $h(\xi)=0$ and $h(\chi)=0$. Thus
$\xi\circ \chi [-1]=0$ thanks to the fact that $\Ii^2=0$. Therefore
there exists $\lambda:X[-1]\to V$ such that $\mu\circ
\lambda=\chi[-1]$, in other words one has the following
commutative diagram
$$
\xymatrix{&&X[-1]\ar[dl]_\lambda\ar[d]^{\chi[-1]}&\\
U\ar[r]^{\alpha} & V\ar[r]^{\mu}& Y\ar[r]^{\xi}& U[1]}
$$
We claim that one can always find $\lambda$ with property $h(\lambda)=0$.
Indeed, for a given $\lambda$ with $\mu\circ \lambda=\chi[-1]$ one obtains the
following diagram after applying $h$:
$$
\xymatrix{&&&x[-1]\ar[dl]_{\hat{\lambda}}\ar[d]^{0}&\\
0\ar[r] &u\ar[r]^{\hat{\alpha}} & v\ar[r]^{\hat{\mu}}& y\ar[r]& 0}
$$
Thus $\hat{\lambda}=\hat{\alpha}\circ \hat{\phi}$, for some
$\phi:X[-1]\to U$. Now it is clear that $\lambda'=\lambda-
\alpha\circ \phi$ has the expected properties $h(\lambda')=0$ and
$\mu\circ \lambda'=\chi[-1]$, and the claim is proved.

One can use the morphisms axiom to conclude that there exists a
commutative diagram
$$\xymatrix{X[-1]\ar[r]^\lambda \ar[d]^{\id}& V \ar[r] \ar[d]^\mu &Z\ar[r]\ar@{-->}[d]
& X\ar[d]^\id\\
X[-1]\ar[r]^{\chi[-1]}& Y\ar[r]^\nu  &W\ar[r]^\gamma & X}
$$
Since  $h(\lambda)=0$, by applying $h$ one obtains the following commutative
diagram
$$\xymatrix{0\ar[r]&v\ar[r]\ar[d]_{\hat{\mu}} &z\ar[r]\ar[d] &x\ar[r]\ar[d]_\id&0\\
0\ar[r]&y\ar[r]^{\hat{\nu}}&w\ar[r]^{\hat{\gamma}}&x\ar[r]&0}
$$
which shows that one has a commutative diagram with exact rows
$$\xymatrix{&0\ar[r]&v\ar[r]\ar[d]_{\id} &z\ar[r]\ar[d] &x\ar[r]\ar[d]_\id&0\\
0\ar[r]&u\ar[r]^{\hat{\alpha}}&v\ar[r]^{\hat{\beta}}&
 w\ar[r]^{\hat{\gamma}}&x\ar[r]&0}
$$
Thus $\sA$ is  hereditary.

\smallskip
\noindent {\it Forth step. Idempotents split in $\sT$}.  Let ${\sf
Idem}(\sT)$ be the idempotent completion of  $\sT$ (see
\cite{arefree} or \cite{ictc}). We have to show that the canonical
functor $\sT\to {\sf Idem}(\sT)$ is an equivalence of categories.
One can summarize the previous steps saying that the category
$\sT$ is a linear extension of $\sA$ by the bifunctor
$(X,Y)\mapsto \ext^1_{\sA}(h(X)[1],h(Y))$ in the sense of Baues
and Wirsching \cite{bw}. Now one can use Proposition 3.2 of
\cite{arefree} to conclude that $\sT\to {\sf Idem}(\sT)$ is indeed
an equivalence of categories.

An alternative proof can be done using the result of  \cite{ictc}
and Corollary \ref{co2} below which uses only the first three
steps. Indeed, by \cite{ictc}, the category $\sT'={\sf Idem}(\sT)$
carries a natural triangulated structure. Since $\sA$ is an
abelian category, all idempotents in $\sA$ split and it follows
from the universal property of the idempotent completion that the
functor $h$ has a unique extension $\sT'\to \sA$, which
is  denoted by $h'$.  We claim that the functor $h'$ reflects isomorphisms. Indeed, if
$X'$ is an object in $\sT$ such that $h'(X')=0$, then there exits an object $Y'$ such that
$Z=X'\oplus Y'$ lies in $\sT$.  Let $e:Z\to Z$ be given by $e(x,y)=(0,y)$. Then $h(Z)=h'(Y')$ and therefore $h(e)$ is an isomorphism. By our assumption on $h$ it follows that
$e$ is an isomorphism and hence $X'=0$. It is clear that $h'$ is full and realizes all objects of $\sA$. Hence  the conditions of
Corollary \ref{co2} below hold and therefore $\sT\to {\sf
Idem}(\sT)$ is an equivalence of categories. 
\end{proof}

\begin{Le}\label{mzmrelizacia} Let $h:\sT\to \sA$ be a homology theory. Assume $h$ reflects
isomorphisms and is full. Suppose there is given a morphism $\alpha:U\to V$, an
object $W$ in $\sT$ and a short exact sequence
$$0\to u \buildrel {\hat{\alpha}} \over \to  v \buildrel {\hat{\beta}} \over \to  w\to 0 
$$
in $\sA$, where as usual $u=h(U), v=h(V), w=h(W)$ and
$\hat{\alpha}=h(\alpha)$. Then there exists a distinguished triangle
$$U \buildrel {\alpha} \over \to  V \buildrel {\beta} \over \to  W \buildrel {\gamma} \over \to  U[1]
$$
such that $h(\beta)=\hat{\beta}$. The dual statement is also true: Suppose
there is given a morphism $\beta:V\to W$, an object $U$ in $\sT$ and a short
exact sequence
$$0\to u \buildrel {\hat{\alpha}} \over \to  v \buildrel {\hat{\beta}} \over \to  w\to 0
$$
in $\sA$, where  $\hat{\beta}=h(\beta)$. Then there exists a
distinguished triangle
$$U \buildrel {\alpha} \over \to  V \buildrel {\beta} \over \to  W \buildrel {\gamma} \over \to  U[1]
$$
such that $h(\alpha)=\hat{\alpha}$.
\end{Le}

\medskip

\begin{proof} Take any distinguished triangle containing $\alpha$,
$$U \buildrel {\alpha} \over \to  V \buildrel {\eta} \over \to  Z \buildrel {\epsilon} \over \to  U[1].
$$
Apply $h$ to get a short exact sequence
$$0\to u \buildrel {\hat{\alpha}} \over \to  v \buildrel {\hat{\eta}} \over \to  z\to 0.
$$
Then we get the following commutative diagram
$$\xymatrix{0\ar[r]& u\ar[r]^{\hat{\alpha}} \ar[d]^\id &
v\ar[r]^{\hat{\eta}}\ar[d]^\id &z\ar[d]^{\hat{\delta}}\ar[r]&0\\
0\ar[r]& u\ar[r]^{\hat{\alpha}} & v\ar[r]^{\hat{\beta}}&w\ar[r] &0}$$
with $\hat{\delta}$ an isomorphism. By assumption one can realize
$\hat{\delta}$ to obtain an isomorphism $\delta:Z\to W$,
$h(\delta)=\hat{\delta}$. Then we have an isomorphism of triangles
$$\xymatrix{U\ar[r]^{\alpha}\ar[d]^\id  &
V\ar[r]^{\eta} \ar[d]^\id &Z\ar[r]^\epsilon \ar[d]^\delta&U[1]\ar[d]^\id\\
U\ar[r]^{\alpha} & V\ar[r]^{\beta}&W\ar[r]^\gamma &U[1]}
$$
where $\beta=\delta\eta$ and $\gamma=\epsilon\circ \delta^{-1}$.
It follows that the triangle
$$U \buildrel {\alpha} \over \to  V \buildrel {\beta} \over \to W \buildrel {\gamma} \over \to U[1]
$$
is also a distinguished triangle. Thus the first statement is proved. The  dual
argument gives the second result. 
\end{proof}

\begin{Co}\label{co2} Let $j:\sT\to \sT'$ be a triangulated functor between triangulated categories.
Assume $h':\sT'\to\sA$ is a homological functor satisfying the
conditions i), ii) and iv) of Theorem \ref{teorema}. If the homology
functor $h=h'\circ j:\sT\to \sA$ also satisfies the same
conditions then $j$ is an equivalence of categories.
\end{Co}

\medskip
\begin{proof} First observe that the functor $j$ is full and faithful
because for any pair of objects $X,Y\in \sT$ both abelian groups
$\sT(X,Y)$ and $\sT'(jX,jY)$ are part of the equivalent extensions
of $\hom_{\sA}(h(X),h(Y))$ by $\ext^1_{\sA}(h(X)[1],h(Y))$. If now
$X'$ is an object in $\sT'$ then there is an object $X$ in $\sT$
and an isomorphism $\hat{\alpha}:h(X)\to h'(X')$ in $\sA$. But
$h(X)=h'(j(X))$ and $h'$ is full so $\hat{\alpha}=h'(\alpha)$ for
a morphism $\alpha:jX\to X'$, which is an isomorphism because $h'$
reflects isomorphisms.
\end{proof}

\bigskip
\noindent
{\bf ACKNOWLEDGEMENTS.}
The first author was supported by the University of Bielefeld and C.N.R.S.  He also acknowledge  discussions with Vincent Franjou, Bernhard Keller, Claus Michael Ringel and Stefan Schwede.


\begin{thebibliography}{99}

\bibitem{ictc} {\sc P. Balmer} and {\sc M. Schlichting}.
Idempotent completion of triangulated categories.
J. of Algebra,   236  (2001) , 819--834.

\bibitem{bw} {\sc H.-J. Baues} and {\sc G. Wirsching}.
Cohomology of small categories.
J. Pure  Appl. Algebra  38 (1985), 187--211.

\bibitem{phantom} {\sc J.D. Christensen and N. P. Strickland.}
Phantom maps and homology theories.  Topology  37  (1998),  no. 2, 339--364.

\bibitem{GM} {\sc S. I. Gelfand and Y. I. Manin}.
Methods of homological algebra. Second edition. Springer
Monographs in Mathematics. Springer-Verlag, Berlin, 2003. xx+372
pp.

\bibitem{arefree} {\sc T. Pirashvili}.
Projectives are free for nilpotent algebraic theories.
Algebraic $K$-theory and its applications (Trieste, 1997),
589--599, World Sci. Publishing, River Edge, NJ, 1999.

\bibitem{puppe} {\sc D. Puppe}. On the structure of stable homotopy theory. Colloquium
on algebraic topology. Aarhus Universitet Matematisk Institut
(1962), 65--71.

\bibitem{verdier} {\sc J. L. Verdier}. Des Cat\'egories deriv\'ees.
des cat\'egories ab\'eliennes. Ast\'erisque  v. 239 (1996). 253
pp.

\end{thebibliography}
\end{document}